\documentclass[letterpaper,twocolumn,10pt]{article}
\usepackage{usenix}

\usepackage[available,functional,reproduced]{usenixbadges}

\usepackage{algorithm}
\usepackage[noend]{algpseudocode}
\usepackage{amssymb}
\usepackage{amsmath}
\usepackage{booktabs}
\usepackage{listings}
\usepackage{graphicx}
\usepackage{hyperref}
\usepackage{caption}
\usepackage{subcaption}
\usepackage{rotating}
\usepackage{enumitem}
\usepackage{nopageno}
\setitemize{noitemsep,topsep=1.5pt}
\setenumerate{noitemsep,topsep=1.5pt}

\begin{document}

\title{Revealing Floating-Point Accumulation Orders in Software/Hardware Implementations}

\author{
{\rm Peichen Xie}\\
Microsoft Research
\and
{\rm Yanjie Gao}\\
Microsoft Research
\and
{\rm Yang Wang}\\
Microsoft Research
\and
{\rm Jilong Xue}\\
Microsoft Research
} 

\thispagestyle{plain}
\pagestyle{plain}

\maketitle
\begin{abstract}
Accumulation-based operations, such as summation and matrix multiplication, are fundamental to numerous computational domains. However, their accumulation orders are often undocumented in existing software and hardware implementations, making it difficult for developers to ensure consistent results across systems. 
To address this issue, we introduce FPRev, a diagnostic tool designed to reveal the accumulation order in the software and hardware implementations through numerical testing. With FPRev, developers can identify and compare accumulation orders, enabling developers to create reproducible software and verify implementation equivalence.

FPRev is a testing-based tool that non-intrusively reveals the accumulation order by analyzing the outputs of the tested implementation for distinct specially designed inputs. Employing FPRev, we showcase the accumulation orders of popular libraries (such as NumPy and PyTorch) on CPUs and GPUs (including GPUs with specialized matrix accelerators such as Tensor Cores). We also validate the efficiency of FPRev through extensive experiments. FPRev exhibits a lower time complexity compared to the basic solution. FPRev is open-sourced at \url{https://github.com/peichenxie/FPRev}.

\end{abstract}

\section{Introduction}

Today, floating-point computations are ubiquitous, with accumulation-based operations (AccumOps) such as summation, dot products, matrix-vector multiplications, and matrix multiplications playing fundamental roles in various domains. However, no general specification dictates the accumulation orders of AccumOps. Without well-defined requirements, AccumOps are implemented differently across software and hardware, leading to inconsistencies due to the non-associativity of floating-point addition \cite{villa_effects_2009}. For example, the half-precision (float16) sum of 0.5, 512, and 512.5 depends on the accumulation order: $(0.5+512)+512.5=1025$, while $0.5+(512+512.5)=1024$. Consequently, varying AccumOp implementations yield different results, complicating reproducibility in software development.

Numerical reproducibility is critical in scientific computing \cite{he_using_2001,taufer_improving_2010,pouchard_computational_2019}, high-performance computing \cite{villa_effects_2009}, database management system \cite{muller_reproducible_2018}, deep learning \cite{shanmugavelu_impacts_2024,chen_towards_2022}, etc. Particularly, software without verified numerical reproducibility is deemed risky or disqualified when applied to safety-critical or rigorous scenarios like aerospace or banking, where even minor inconsistencies in data are unacceptable. Unfortunately, with the rapid evolution of heterogeneous hardware and the fast iteration of diverse software stacks, reproducibility of AccumOps has become increasingly challenging. Existing implementations rarely disclose their accumulation orders, hindering reproducible AccumOp development.



We propose FPRev, a diagnostic tool to help developers identify how AccumOps are implemented in software and hardware. FPRev reveals the accumulation order of an AccumOp implementation (AccumImpl) through numerical testing. This enables developers to reproduce an AccumImpl on a new system by using the revealed accumulation order as a specification and verify equivalence between two implementations by comparing their accumulation orders.

\begin{figure*}[t]
    \centering
    \includegraphics[width=\linewidth]{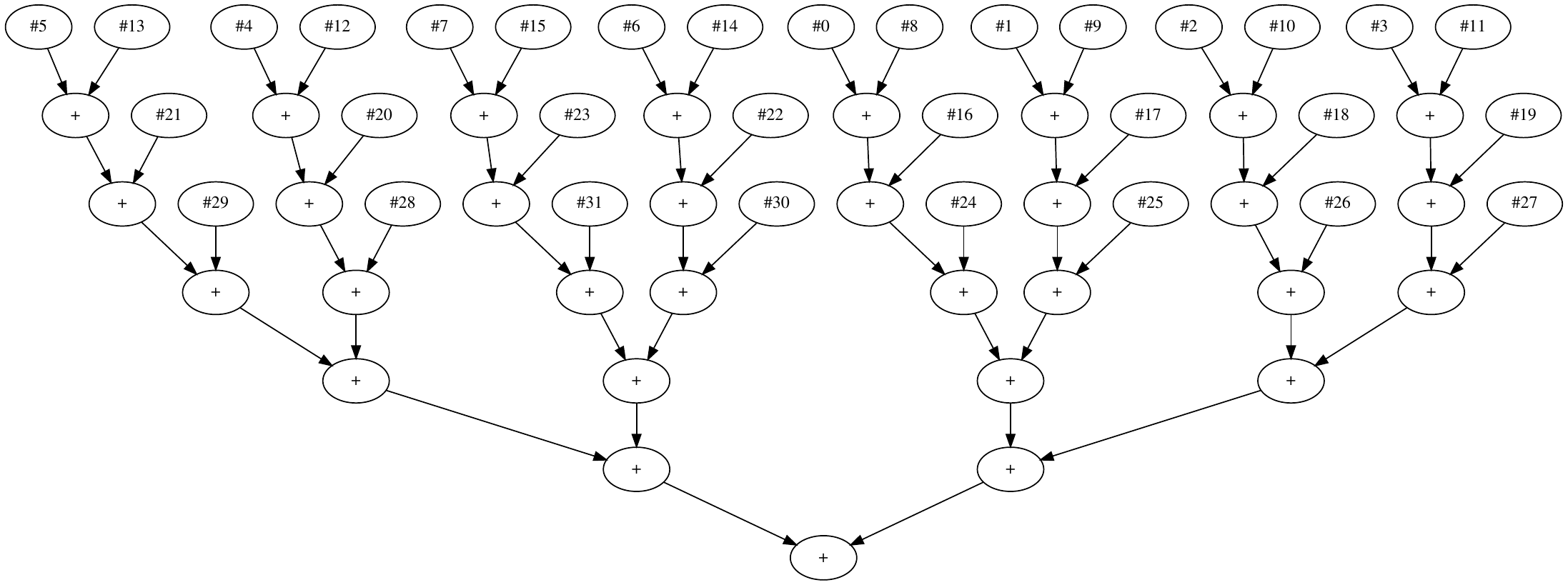}
    \caption{Visualizing the accumulation order of Numpy's summation function for $n=32$ single-precision numbers with a summation tree. The numbers on the leaf nodes denote the indexes in the input.}
    \label{fig:numpy_sum_32}
\end{figure*}

As a case study, we use FPRev to analyze popular numerical libraries on diverse hardware, uncovering their undocumented and undisclosed accumulation orders. On different CPUs, we apply FPRev to the NumPy library \cite{harris_array_2020}. On different GPUs, we apply FPRev to the PyTorch library \cite{paszke_pytorch_2019}. The results indicate that NumPy's summation functions are implemented equivalently across CPUs, and the same holds for PyTorch's summation functions across GPUs. However, other AccumOps relying on BLAS (Basic Linear Algebra Subprograms) backends like Intel MKL \cite{mkl}, OpenBLAS \cite{OpenBLAS}, and NVIDIA cuBLAS \cite{cublas} exhibit non-reproducible behavior.

FPRev also visualizes the order with the summation tree, i.e., a full binary tree representing how an AccumImpl performs summation, to guide develeopment. For example, Figure \ref{fig:numpy_sum_32} illustrates NumPy's summation of 32 single-precision (float32) numbers. It divides the 32 numbers into 8 ways, accumulates the summands with a stride of 8 on each way, and sums up the 8 ways together using pairwise summation. This 8-way accumulation order is friendly to CPU's SIMD instructions. With this information, it is easy to replicate NumPy's numerical behavior in a new implementation. 

\textbf{Design overview.} Determining the accumulation order of an AccumImpl is a challenging task. Static methods, such as analyzing the source code, are cumbersome and inapplicable to black-box implementations and compiler optimization. Dynamic methods, like scrutinizing the runtime traces, lack an automatic tool to analyze the traces. In addition, many software or hardware implementations are parallel, making the analysis more challenging.

We address these challenges through non-intrusive testing. Recall the example where $(0.5+512)+512.5=1025$ and $0.5+(512+512.5)=1024$ in half precision. Different accumulation orders yield distinct results, making it possible to deduce the order from numerical outputs. However, the number of all possible accumulation orders is exponential, making the time complexity of the naive brute-force solution (NaiveSol) impractical.

To achieve practical time complexity, we propose a basic solution called BasicFPRev that uses specially designed inputs to facilitate the distinguishing process. We take the summation function for example. First, BasicFPRev set all $n$ summands to 1.0. Then, two of the summands are replaced by a very large number (denoted by $M$) and its negative ($-M$), where $M$ satisfies $(n-2)+M=M$. The summation output corresponds to an integer between 0 and $n-2$, depending on when $M$ or $-M$ cancel each other during the accumulation. Specifically, when $\pm M$ is added to other numbers, it remains $\pm M$; when $M$ is added to $-M$, it results in 0; after that, the remaining summands are accumulated without rounding errors because they are all 1.0. Therefore, the output equals the number of summands accumulated after $M + (-M)$.

BasicFPRev leverages this information to construct the summation tree. Using $i$ and $j$ to denote the indexes of $M$ or $-M$ in the input, we note that the operation $M+(-M)$ corresponds to the lowest common ancestor (LCA) of node \#$i$ and \#$j$ in the summation tree, and the number of leaf nodes under the LCA equals $n$ minus the summation output. Based on this finding, BasicFPRev enumerates $i$ and $j$, collects the output for the corresponding input, and infers the size of subtree rooted at the LCA of node \#$i$ and \#$j$. BasicFPRev then constructs the summation tree bottom-up, starting with subtrees of two leaf nodes and progressively building larger subtrees until the entire tree is generated.


Based on BasicFPRev, FPRev further reduces time complexity by eliminating redundancy, and adds support for matrix accelerators such as Tensor Cores on recent NVIDIA GPUs \cite{markidis_nvidia_2018}. Matrix accelerators are specialized hardware units on GPUs for high-performance matrix multiplication, but they perform non-standard multi-term fused summations \cite{fasi_numerical_2021}. FPRev models their accumulation orders using a multiway tree, where a node with multiple children represents a multi-term fused summation for a group of summands. The summands are aligned and truncated before they are accumulated, as if they are added in finite-precision fixed-point arithmetic.

FPRev has a time complexity of $\Omega(nt(n))$ and $O(n^2t(n))$, where $t(n)$ is the time complexity of the tested AccumImpl. This shows a significant improvement over the $O(4^n/n^{3/2}\cdot t(n))$ complexity of NaiveSol and the $\Theta(n^2t(n))$ complexity of BasicFPRev. Experimental results confirm FPRev’s efficiency and scalability across diverse AccumOp implementations on three CPUs and three GPUs with distinct architectures. 

In summary, the contributions of this paper include:
\begin{enumerate}

\item Design and development of FPRev: we introduce FPRev, a diagnostic tool that non-intrusively reveals the accumulation order of accumulation-based operations implemented in different software and hardware, enabling developers to verify equivalence and maintain numerical reproducibility between implementations.

\item Empirical analysis of popular implementations: We demonstrate FPRev's capabilities by analyzing accumulation orders in popular libraries (e.g., NumPy and PyTorch) across CPUs and GPUs, providing reproducibility insights in backend implementations.

\item Algorithmic innovation: we describe the algorithm of FPRev for revealing accumulation orders and constructing summation trees, refine the algorithm to reduce time complexity, and extend it to handle modern GPU matrix accelerators (e.g., NVIDIA Tensor Cores), modeling their multi-term fused summation using multiway trees.

\item Performance evaluation: we evaluate FPRev's efficiency through comprehensive experiments, test diverse AccumOps implementations on various CPUs and GPUs, and demonstrate significant performance improvements over naive and basic solutions.

\end{enumerate}

\section{Related work}

\subsection{AccumOp implementations}

\subsubsection{On canonical CPUs and GPUs}

Accumulation-based operations (AccumOps) are implemented diversely on modern systems. On most CPUs and GPUs, implementations use standard IEEE-754 addition or fused multiply-add (FMA) arithmetic \cite{noauthor_ieee_2019} to accumulate floating-point numbers. However, they may perform accumulations in different orders without explicitly disclosing those orders.

First, there is diverse numerical software, including BLAS libraries such as Intel MKL \cite{mkl} and NVIDIA cuBLAS \cite{cublas}, Python libraries such as NumPy \cite{harris_array_2020} and PyTorch \cite{paszke_pytorch_2019}, and domain-specific compilers such as Numba \cite{lam_numba_2015} and Triton \cite{tillet_triton_2019}. These libraries are developed without a unified specification, making it difficult to guarantee consistent accumulation orders.

Second, the same software may behave differently across different hardware. Different CPUs and GPUs vary in architecture, number of cores, SIMD width, cache size, etc. Consequently, for performance optimization, software may adjust the accumulation order based on the specific hardware characteristic. Specifically, library developers implement various techniques (e.g., different configurations of loop unrolling, block partitioning, cache optimization, and vectorization) for performance tuning, resulting in different accumulation orders. Additionally, auto-tuners (e.g., Triton \cite{tillet_triton_2019} and TVM \cite{chen_tvm_2018}) are often used to search for optimal configurations, given the complexity of performance factors such as instruction pipelining and dynamic frequency scaling.

Although order-independent algorithms have been proposed \cite{demmel_fast_2013,demmel_parallel_2015,collange_numerical_2015}, which ensure consistent results regardless of the accumulation order, they are highly inefficient and thus rarely used in industry.

Our tool FPRev supports the AccumOp implementations on canonical CPUs and GPUs and can reveal their undisclosed accumulation orders.

\subsubsection{On matrix accelerators}

Matrix accelerators \cite{markidis_nvidia_2018,raihan_modeling_2019} are specialized hardware components in modern GPUs designed for high-performance matrix multiplication. 
Developers can implement matrix multiplications using the APIs of matrix accelerators \cite{ptx}. However, the numerical behavior and accumulation order of the APIs are undocumented and inconsistent across different GPUs.

FPRev supports the AccumOp implementations based on matrix accelerators and can reveal their undisclosed accumulation orders.

\subsection{Revelation of numerical behaviors}

FPRev achieves non-intrusive revelation of accumulation orders through numerical testing. Prior works \cite{fasi_numerical_2021,li_fttn_2024} have also employed numerical testing to study the numerical behavior of matrix accelerators. They design ``corner cases'' as test inputs and analyze the numerical behavior based on the outputs. They find that for float64 on NVIDIA Tensor Cores and AMD Matrix Cores, the matrix multiplication instruction is based on a chain of standard FMA arithmetic. In contrast, other instructions use a non-standard arithmetic where multiple terms (the exact number depends on the hardware) are accumulated after alignment and truncation.

FPRev is a general tool that applies to AccumOp implementations, including those based on matrix accelerators, while prior works focus exclusively on specific hardware.

\subsection{Numerical reproducibility engineering}

The inconsistent AccumOp implementations pose significant issues in numerical reproducibility \cite{villa_effects_2009,bailey_facilitating_2016}. To help developers debug the issues, several testing-based tools have been proposed. For example, Varity \cite{laguna_varity_2020} uses randomized testing to verify equivalence between implementations. Tools like pLiner \cite{guo_pliner_2020} and its follow-up \cite{miao_expression_2023} employ differential testing to pinpoint non-reproducible parts of a program. In contrast, FPRev uses a deterministic testing method to identify the accumulation order of AccumOps.

Early works \cite{he_using_2001,taufer_improving_2010} have emphasized the importance of reproducible AccumOps for ensuring numerical stability but lack practical solutions. A preliminary approach \cite{arteaga_designing_2014} adopts the aforementioned order-independent algorithm \cite{demmel_fast_2013} to ensure reproducibility, but suffers from its inefficiency. In contrast, FPRev offers a more practical solution: replicating the accumulation order of existing efficient implementations.

\section{Problem statement}
\label{sec:problem}

\subsection{Motivation}

Accumulation-based operations (AccumOps) are fundamental in floating-point computing, but most implementations do not specify their accumulation orders. This lack of transparency motivates us to design a tool for revealing the accumulation order.

An example application of the tool is in developing reproducible AccumOps, which are key to software reproducibility and service consistency \cite{he_using_2001,taufer_improving_2010,openai_reproducible_outputs}. Developers must ensure that the accumulation order remains consistent across systems to maintain the reproducibility of AccumOps. If developers can determine the accumulation order, they can use it as a specification to guide their development process. In addition, when porting software to a new system, developers need a rigorous way to verify the equivalence of AccumOps between two systems. This can be achieved by comparing the accumulation orders of the AccumOps implemented on two systems. 

\subsection{Problem definition}

For an AccumOp implementation (AccumImpl), we aim to design a diagnostic tool to reveal its accumulation order. For brevity, we focus on the summation function in the following discussion, since other AccumOps can be abstracted as calls to the summation function with the intermediate results as inputs. For example, dot product $\mathbf{x\cdot y}$ can be treated as $\sum_{i=0}^{n-1} x_iy_i$. Thus, solutions for the summation function can be naturally applied to other AccumOps.

We formulate the summation operation as follows. The floating-point addition is performed $n-1$ times in a predetermined order to calculate the sum of $n$ floating-point numbers. We assume that the accumulation order is unknown but is uniquely determined by the given implementation on specific hardware. Therefore, randomized implementations and those where the order depends on the values of the summands are out of scope\footnote{While randomized or input-dependent summation algorithms are possible, they tend to be inefficient. Similarly, algorithms using AtomicAdd are affected by thread scheduling, but we did not find popular libraries implementing such algorithms.}.

The problem is how to reveal the accumulation order in an summation implementation, denoted by \textsc{SumImpl}. Specifically, the input of our revelation algorithm is \textsc{SumImpl} and the number of summands $n$. The output of the algorithm is the accumulation order of \textsc{SumImpl}. 

Strictly speaking, the accumulation order is represented by a computational graph called the summation tree, which is a rooted full binary tree with $n$ leaf nodes and $n-1$ inner nodes. Each addition operation corresponds to an inner node, which represents the sum of this operation. The two children of the node represent the two summands of this operation. For example, Figure \ref{fig:numpy_sum_32} depicts the accumulation order of Numpy's summation function for $n=32$ single-precision numbers.



\subsection{Inefficiency of the naive solution}
\label{sec:naive}
We now introduce a naive solution (NaiveSol) to the problem, which is based on brute-force search. We design a recursive algorithm to enumerate every possible accumulation orders. For each order, we verify its correctness through randomized testing. Specifically, we generate multiple random inputs, compute the sums in the current order, and compare the results with those from  \textsc{SumImpl}. If the results match, we accept the order.


The time complexity of NaiveSol is $O(4^n/n^{3/2}\cdot t(n))$, as the number of all possible orders is the $(n-1)$-th Catalan number $C_{n-1}=\frac{(2n-2)!}{n!(n-1)!}=O(4^n/n^{3/2})$. Here, $t(n)$ represents the time complexity of \textsc{SumImpl}. In addition to being inefficient, NaiveSol is not fully reliable because different orders can produce the same output for certain inputs. Although the probability is low and reliability can be improved by increasing the number of test inputs, a deterministic solution with full reliability is preferable, as we will achieve next.


\section{Basic polynomial-time solution}
\label{sec:method}

The exponential complexity of the naive solution is highly impractical. To address the issue, we present our basic solution called BasicFPRev for revealing the accumulation order, which reduces the time complexity to polynomial. We design an algorithm to determine the accumulation order from the numerical results of the tested summation implementation  (\textsc{SumImpl}) for specially designed testing inputs. The following parts detail the three steps of the algorithm.

\subsection{Step 1: designing testing inputs}

To facilitate distinguishing the accumulation order, we leverage the swamping phenomenon of floating-point addition \cite{higham_accuracy_1993}. When two floating-point numbers differing by many orders of magnitude are added, the smaller number is swamped and makes no contribution to the sum. For example, $2^{24}+1$ equals $2^{24}$ in single-precision (float32) arithmetic. 

To induce and utilize this phenomenon, we construct various ``masked all-one arrays" as testing inputs. Specifically, let $n$ denote the number of summands, and let \textsc{SumImpl} represent an summation implementation with a predetermined but unknown accumulation order. Let $M$ be a very large floating-point number that readily induces the swamping phenomenon. For example, we set $M=2^{127}$ for float32 or $M=2^{1023}$ for float64. Then, we define a masked all-one array $A^{i,j}$ as
$A^{i,j} = (A_0^{i,j}, A_1^{i,j}, ..., A_{n-1}^{i,j})$
such that
$$ A_k^{i,j}=\begin{cases}
    M   & \quad \text{if } k=i\\
    -M  & \quad \text{if } k=j\\
    1.0 & \quad \text{otherwise}
\end{cases}$$
where $i$ and $j$ denote the indexes of $M$ and $-M$ in the array. In $A^{i,j}$, there exist exactly one $M$ and one $-M$, with all other elements being $1.0$.

We use $\pm M$ as masks. Specifically, $M+\sigma = M$ and $-M+\sigma = -M$ hold for $0\le \sigma \le n-2$ in floating-point arithmetic, if $n\ll M$. Therefore, in \Call{SumImpl}{$A^{i,j}$}, adding any summand or intermediate sum (except $M$ and $-M$ themselves) to $\pm M$ yields $\pm M$. In other words, $M$ and $-M$ can mask the summands or intermediate sums added to them.

As a result, the output of \Call{SumImpl}{$A^{i,j}$} depends on the accumulation order and we can distinguish the accumulation order from the output. For example, given $n=3$ and $A^{0,1}=(M, -M, 1)$, sequential summation $M+(-M)+1$ corresponds to 1, stride summation $M+1+(-M)$ corresponds to 0, and reverse summation $1+(-M)+M$ corresponds to 0. If the output equals 1, then we can infer the accumulation order is sequential summation. If the output equals 0, we can determine the exact accumulation order by further testing with $A^{0,2}$ as the input.


\subsection{Step 2: analyzing the accumulation order from the outputs}

To analyze the accumulation order, we call \textsc{SumImpl} with $n(n-1)/2$ inputs, i.e., $A^{i,j}$ for $0\le i<j<n$. Each output reveals information about the accumulation order. Specifically, since $\pm M$ mask the summands or intermediate sums added to them, these numbers make no contribution to the sum. In contrast, only those summands not masked by $\pm M$ contribute to the sum. Therefore, the output equals the sum of these summands. Since each of the summands equals $1.0$, the output equals the number of the summands not masked by $\pm M$:
$$n_\mathrm{not\ masked}^{i,j} = \Call{SumImpl}{A^{i,j}}.$$
Then, we can also obtain the number of the summands masked by $\pm M$ by calculating $n_\mathrm{masked}^{i,j}=n-2-n_\mathrm{not\ masked}^{i,j}$.

How does this information relate to the order, or specifically, the summation tree? Recall that $i$ and $j$ denote the positions of the masks, represented by node  \#$i$ and \#$j$ in the summation tree. \textbf{We note that the neutralization of the two masks (i.e., the addition operation $M+(-M)=0$) corresponds to the lowest common ancestor (LCA) of node  \#$i$ and \#$j$ }. Then, observing the subtree rooted at the LCA, we find that all the summands masked by $\pm M$ are in the subtree, and all the summands not masked by $\pm M$ are out of the subtree. Therefore, the number of leaf nodes in the subtree (representing the size of the subtree) equals $n-n_\mathrm{not\ masked}^{i,j}$, denoted by
$$l^{i,j}=n-n_\mathrm{not\ masked}^{i,j}=n-\Call{SumImpl}{A^{i,j}}.$$
For brevity, we use $l^{i,j}$ to denote ``the number of leaf nodes in the subtree rooted at the LCA of node  \#$i$ and \#$j$ '' in the rest of the paper.

\begin{figure}[h]
    \centering
    \includegraphics[width=0.85\columnwidth]{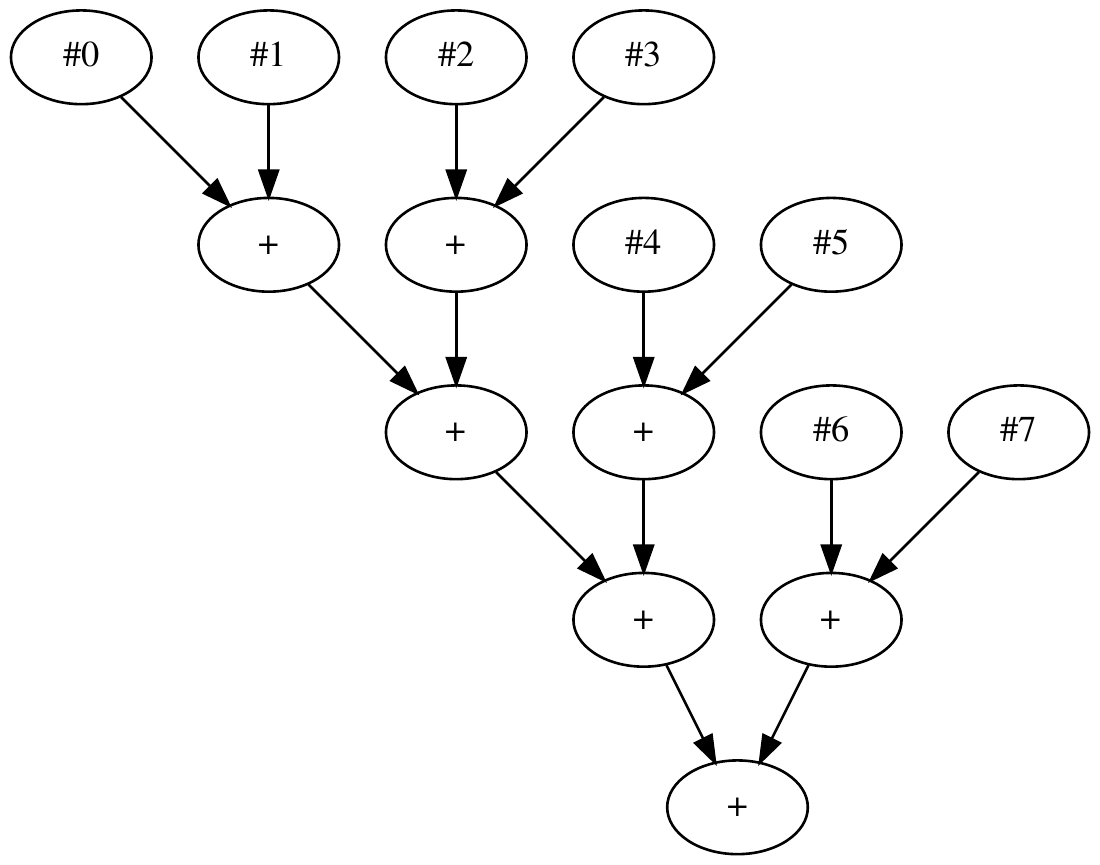}
    \caption{The summation tree of Algorithm \ref{alg:sum8}. The numbers on the leaf nodes denote the indexes in the input.}
    \label{fig:sum8}
\end{figure}

Take Algorithm \ref{alg:sum8} as an example \textsc{SumImpl}, whose accumulation order is depicted in Figure \ref{fig:sum8}.\footnote{We compile it with different compiler versions and optimization flags, and note that they do not change the accumulation order, even if they might in theory. We suspect such changes can be found on more complex math expressions instead of summations.} If $i=2$ and $j=4$, then the array $A^{2,4}$ is $(1,1,M,1,-M,1,1,1)$. Computing the sum of $A^{2,4}$ with the example \textsc{SumImpl}, the 3rd summand and the intermediate sum of the 0th and 1st summands are masked by $M$ (the 2nd summand); the 5th summand is masked by $-M$ (the 4th summand). Therefore, in total, $n_\mathrm{masked}^{2,4}=4$. In contrast, the 6th and 7th summands and their intermediate sum are not added to $M$ or $-M$, so $n_\mathrm{not\ masked}^{2,4}=\Call{SumImpl}{A^{2,4}}=2$.

\begin{algorithm}[h]
\caption{An example summation implementation.}
\label{alg:sum8}
\begin{lstlisting}{language=C++} 
float sum = 0;
for (int i=0; i<8; i+=2)
    sum += a[i] + a[i+1];
\end{lstlisting}
\end{algorithm}

The LCA of node \#2 and \#4 is the grandparent node of node \#4, as shown in \ref{fig:sum8}. It corresponds to the neutralization of the 2nd summand $M$ and the 4th summand $-M$, i.e., $M+(-M)=0$. Within the subtree rooted there, there are node \#0, \#1, \#2, \#3, \#4, and \#5 (6 leaf nodes in total), corresponding to the two masks and the summands masked by them. In contrast, node \#6 and \#7 (2 leaf nodes in total), which correspond to the summands not masked, are out of the subtree. Therefore, $l^{2,4}=8-2=6$. Table \ref{tab:sum8} shows more examples of the output of \Call{SumImpl}{$A^{i,j}$} and $l^{i,j}$ for Algorithm \ref{alg:sum8}.

\begin{table}[ht]
\caption{The order-related information $l^{i,j}$ inferred from the outputs of Algorithm \ref{alg:sum8} with different masked all-one arrays $A^{i,j}$ as inputs.}
\label{tab:sum8}
\centering
\small
\begin{tabular}{@{}ccccc@{}}
\toprule
$i$ & $j$ & Input: $A^{i,j}$& Output& $l^{i,j}$ \\ \midrule
0 & 1 & $(M,-M,1,1,1,1,1,1)$ & 6 & 2 \\
0 & 2 & $(M,1,-M,1,1,1,1,1)$ & 4 & 4 \\
0 & 3 & $(M,1,1,-M,1,1,1,1)$ & 4 & 4 \\
0 & 4 & $(M,1,1,1,-M,1,1,1)$ & 2 & 6 \\
0 & 5 & $(M,1,1,1,1,-M,1,1)$ & 2 & 6 \\
0 & 6 & $(M,1,1,1,1,1,-M,1)$ & 0 & 8 \\
0 & 7 & $(M,1,1,1,1,1,1,-M)$ & 0 & 8 \\
 &  & $...$ &  &  \\
2 & 3 & $(1,1,M,-M,1,1,1,1)$ & 6 & 2 \\
2 & 4 & $(1,1,M,1,-M,1,1,1)$ & 2 & 6 \\
 &  & $...$ &  &  \\ \bottomrule
\end{tabular}
\end{table}


\subsection{Step 3: generating the summation tree}

With the information $L=\{(l^{i,j}, i, j)\}$ derived from the outputs, where $l^{i,j}$ represents the size of the subtree rooted at the LCA of node  \#$i$ and \#$j$, generating the summation tree from $L$ is a tree algorithm problem.

Our solution employs a bottom-up approach to construct the tree. First, we sort $L$ in ascending order. For each $l^{i,j}$, we locate the root of the existing subtree containing node \#$i$ and the root of the existing subtree containing node \#$j$, and merge them by creating a new parent node for them. By repeating this process, we construct the entire summation tree, progressing from small subtrees to larger ones.

For example, consider the order-related information $L=\{(l^{i,j}, i, j)\}$ shown in Table \ref{tab:sum8}. To generate the summation tree, we start by initializing the tree with eight disjoint nodes labeled 0 to 7. Then, examining the smallest value in $L$, we have $l^{0,1}=2$. This implies that the subtree rooted at the LCA of node \#0 and \#1 should have 2 leaf nodes. Since the summation tree is a full binary tree, node \#0 and \#1 are exactly the two children of the root of this subtree. Therefore, we add a new node to the tree, label it with $n$ plus the label of its left child (i.e., $n+0$ in this example), and add two edges from the two leaf nodes to the new node.

For $l^{2,3}=l^{4,5}=l^{6,7}=2$, similarly, we can construct the subtree containing node \#2 and \#3, the subtree containing node \#4 and \#5, and the subtree containing node \#6 and \#7. Now, we have four subtrees of size 2, where the size of a subtree is represented by the number of leaf nodes in it. 

Next, the smallest unexamined value in $L$ is $l^{0,2}=4$. This implies that the subtree rooted at the LCA of node \#0 and \#2 should have 4 leaf nodes. We note that node \#0 and \#2 are currently in two different subtrees (each has 2 leaf nodes), so we should merge the two subtrees. Therefore, we find the current roots of the subtrees containing node \#0 and \#2 respectively, i.e., node \#$i'$ and \#$j'$ where $i'=n+0$ and $j'=n+2$. Subsequently, we add a new node as the parent of node \#$i'$ and \#$j'$, label it with $n$ plus the label of its left child, i.e., $n+i'=2n+0$, and add two edges from node \#$i'$ and \#$j'$ to the new node.

For $l^{0,3}=4$, we find that node \#0 and \#3 are already in the same subtree with 4 leaf nodes, so we just skip it. The same process is applied to $l^{1,2}=l^{1,3}=4$. Now, we have a subtree with 4 leaf nodes and two subtrees with 2.

The next smallest unexamined value in $L$ is $l^{0,4}=6$, which implies that the subtree rooted at the LCA of node \#0 and \#4 should have 6 leaf nodes. Similarly, node \#0 and \#4 are currently in two different subtrees (one has 4 leaf nodes and the other has 2), so we should merge the two subtrees. Therefore, following the similar process, we find the current roots of the trees containing node \#0 and \#4 respectively, i.e., node \#$i'$ and \#$j'$ where $i'=2n+0$ and $j'=n+4$. Subsequently, we create a new node as their parent (labelled as $n+i'=3n+0$), and add two edges from node \#$i'$ and \#$j'$ to it.

For $l^{0,5}=6$, we find that node \#0 and \#5 are already in the same subtree with 6 leaf nodes, so we skip it. The same process is applied to $i\in \{1,2,3\}$ and $j\in \{4,5\}$. Now, we have a subtree of with 6 leaf nodes and a subtree of with 2.

Finally, in the similar way, the next smallest unexamined value in $L$ is $l^{0,6}=8$, which implies that the subtree rooted at the LCA of node \#0 and \#6 should have8 leaf nodes. Since node \#0 and \#6 are currently in two different subtrees (one has 6 leaf nodes and the other has 2), we should merge the two subtrees. After we add the parent node of the current roots of the two subtrees and add the corresponding edges, the entire summation tree is generated. 

To generalize and formulate the algorithm, we present BasicFPRev in Algorithm \ref{alg:basic}, where the  \textsc{GenerateTree} function encapsulates Step 3. It first initializes $T$ with $n$ disjoint nodes and no edges. Then, for each $l^{i,j}$, it finds the roots of the existing subtrees containing node \#$i$ and \#$j$ . If they are identical, then the $i$-node \#$i$ and \#$j$ are already in the same subtree. Otherwise, it combines them. The FindRoot function can be implemented by the disjoint-set data structure, resulting in an amortized time complexity of $O(\alpha(n))$ \cite{tarjan_worst-case_1984}, where $\alpha(n)$ is the inverse Ackermann function.

\begin{algorithm}[ht]
\caption{BasicFPRev: our basic solution for revealing the accumulation order.}
\label{alg:basic}
\begin{algorithmic}
\Require Implementation \textsc{SumImpl} and the number of summands $n$
\Ensure Summation tree of \textsc{SumImpl}
\Function{BasicFPRev}{\textsc{SumImpl}, $n$}
\State $L\gets \varnothing$ 
\For{$i\gets 0$ to $n-1$ }
    \For{$j\gets i+1$ to $n-1$ }
        \State $A^{i,j} \gets (1,1,...,1)$; $A^{i,j}_i\gets M$; $A^{i,j}_j\gets -M$ 
        \State $l^{i,j} \gets n-\Call{SumImpl}{A^{i,j}}$
        \State $L\gets L\cup \{(l^{i,j},i,j)\}$
    \EndFor
\EndFor
\Function{GenerateTree}{$L$}
    \State $T\gets \varnothing$
    \For{$(l^{i,j},i,j) \in L$ in ascending order}
        \State $i' \gets T.\mathrm{FindRoot}(i)$
        \State $j' \gets T.\mathrm{FindRoot}(j)$
        \If{$i'\neq j'$}   
            \State $k\gets i'+n$ \Comment{assign a new label}\
            \State $T\gets T\cup \{(i',k),(j',k)\}$
        \EndIf
    \EndFor
\State \Return $T$
\EndFunction

\State \Return \Call{GenerateTree}{$L$}
\EndFunction

\end{algorithmic}
\end{algorithm}

\subsection{Time complexity and correctness analysis}
\label{sec:basic-timecor}

Let $t(n)$ denote the time complexity of \textsc{SumImpl}. The computation of $L$ has a time complexity of $\Theta(n^2t(n))$. In \textsc{GenerateTree}, the time complexity of sorting $n(n-1)/2$ elements is $\Theta(n^2\log n^2)=\Theta(n^2\log n)$. Therefore, the time complexity of \textsc{GenerateTree} is $\Theta(n^2\log n + n^2\alpha(n))=\Theta(n^2\log n)$. Thus, the overall time complexity of BasicFPRev is $\Theta(n^2t(n)) + \Theta(n^2\log n)=\Theta(n^2t(n))$.

The correctness of BasicFPRev is inherent in its design  and can be proven as follows. For a given implementation \textsc{SumImpl} and $n$, we use $T$ to denote the real summation tree and define $T'=\Call{BasicFPRev}{\textsc{SumImpl},n}$. Assuming $T\neq T'$, there must exist $i$ and $j$ such that $l^{i,j}_T \neq l^{i,j}_{T'}$. Now we construct $A^{i,j}$ and compute its sum in the order $T$ and $T'$ respectively, resulting in $s$ and $s'$. Then, $s=n-l^{i,j}_T \neq s'=n-l^{i,j}_{T'}$. However, since $s=\Call{SumImpl}{A^{i,j}}$, then $l^{i,j}_{T'}=n-s'\neq n-s=n-\Call{SumImpl}{A^{i,j}}$. This contradicts the statement $l^{i,j} \gets n-\Call{SumImpl}{A^{i,j}}$ in Algorithm \ref{alg:basic}. Therefore, the assumption $T\neq T'$ is false, so $T=T'$.
\section{Algorithm improvement in FPRev}
\label{sec:fprev}

This section introduces the full version of our tool FPRev, which is evolved from our basic solution BasicFPRev detailed in Section \ref{sec:method}. First, we refine the algorithm to reduce its time complexity. Second, based on the refined algorithm, we add support for multi-term fused summation, which is used by matrix accelerators, and finalize FPRev.

\subsection{Reducing time complexity}
\label{sec:optimization}

\subsubsection{Removing redundancy}

By analyzing BasicFPRev, we observe that Algorithm \ref{alg:basic} requires $n(n-1)/2$ different $(l^{i,j},i,j)$ tuples, even though many values of $l^{i,j}$ are identical. However, only $n-1$ new nodes and $2(n-1)$ new edges are constructed. Since computing multiple $l^{i,j}$ by calling \textsc{SumImpl} is the primary source of the method's time complexity, reducing redundancy in $l^{i,j}$ (i.e., the cases corresponding to $i'=j'$ in Algorithm \ref{alg:basic}) can significantly improve efficiency.

To achieve this, we calculate $l^{i,j}$ on demand. Specifically, we do not calculate all $l^{i,j}$ ahead of the tree generation. Instead, we directly start to generate the summation tree, and calculate $l^{i,j}$ when needed. Following the bottom-up idea, we still construct subtrees from leaf to root.

\textbf{Step 1.} We use the set $I=\{0,1,...,n-1\}$ to denote the labels of the leaf nodes for which we are going to construct a summation tree. Let $i$ represent the leaf node with the smallest label in $I$. The sibling node of node \#$i$ is either a leaf node or an inner node. If it is a leaf node, there exists a unique $j$ such that $l^{i,j}=2$. Otherwise, if it is an inner node, then $l^{i,j}>2$ for all $j$ such that $j\neq i$. Therefore, to distinguish the two cases, we need to calculate $l^{i,j}$ for all other $j$s, denoted by the set $L_i=\{l^{i,j}: j\in I-\{i\}\}$. We examine the minimum value in $L_i$, which is denoted by $l = \min(L_i)$.

If $l$ equals $2$, let $j$ be the one that satisfies $l^{i,j}=2$. Then, node \#$j$ is the sibling node of node \#$i$, so we add a new node to the tree, and add two edges from node \#$i$ and \#$j$ to the new node. Now, the currently constructed subtree has 2 leaf nodes.

Otherwise, if $l$ is larger than $2$, the sibling node of node \#$i$ must be an inner node. The subtree rooted at this inner node must have $l-1$ leaf nodes. Let $J_l=\{j:j\in I-\{i\}\land l^{i,j}=l\}$. Then, the number of members of $J_l$ must be $l-1$, and the members of $J_l$ are exactly the leaf nodes of this subtree. This can be proven by contradiction. Now, constructing this subtree is a subproblem for the set $J_l$. Suppose that we have constructed this subtree by a recursive algorithm. We shall add a new node to the tree, and add edges from node \#$i$ and the root node of this subtree to the new node. Now, the currently constructed subtree has $l$ leaf nodes.

Summarizing the two cases, we can treat both cases as the same pattern: finding $J_l=\{j:j\in I-\{i\}\wedge l^{i,j}=l\}$ and solving the subproblem for $J_l$. The first case ($|J_l|=1$) just leads to the stop condition of the recursion ($|I|=1$).

\textbf{Step 2.}  Now we have constructed a subtree with $l$ leaf nodes. Let $r$ be the root of this subtree. Similarly, to find the sibling node of $r$, we examine the minimum value in the rest of $L_i$, which is denoted by $l'$ here. Then, we solve the subproblem for $J_{l'}=\{j:j\in I-\{i\}\wedge l^{i,j}=l'\}$, and get a subtree whose leaf nodes are $J_{l'}$. The root of the subtree, whether a leaf node or an inner node, is the sibling node of $r$. Therefore, we shall add a new node to the tree, and add edges from $r$ and the root node of the subtree to the new node. Now, the currently constructed subtree has $l'$ leaf nodes.

\textbf{Remaining steps.} We repeat the above step until all values in $L_i$ are examined and the entire tree is constructed. We implement this method with a recursive algorithm, as shown in Algorithm \ref{alg:optimized}.

\begin{algorithm}[ht]
\caption{Refinement of BasicFPRev (Algorithm \ref{alg:basic}).}
\label{alg:optimized}
\begin{algorithmic}
\Require Implementation \textsc{SumImpl} and the number of summands $n$
\Ensure Summation tree of \textsc{SumImpl}

\Function{BasicFPRevRefined}{\textsc{SumImpl}, $n$}
    \Function{BuildSubtree}{$I$}
    \State $T\gets \varnothing$
    \If{$|I|=1$} \Comment{stop condition}
        \State \Return $T$
    \EndIf
    \State $i \gets \min(I)$; $L_i\gets \varnothing$
    \For{$j\in I-\{i\}$} \Comment{calculate $l^{i,j}$ on demand}
        \State $A^{i,j} \gets (1,1,...,1)$; $A^{i,j}_i\gets M$; $A^{i,j}_j\gets -M$ 
        \State $l^{i,j}\gets n-\Call{SumImpl}{A^{i,j}}$
        \State $L_i\gets L_i\cup \{l^{i,j}\}$
    \EndFor
    \State $r \gets i$ \Comment{current root of the subtree}
    \For{$l \in L_i$ in ascending order}
        \State $J_l \gets \{j : j\in I-\{i\} \land l^{i,j}=l\}$
        \State $T' \gets \Call{BuildSubtree}{J_l}$ 
        \State $T\gets T\cup T'$
        \State $T\gets T\cup\{(r,r+n), (\mathrm{GetRoot}(T'),r+n)\}$
        \State $r \gets r+n$
    \EndFor
    \State \Return $T$
    \EndFunction
\State \Return \Call{BuildSubtree}{$\{0,1,...,n-1\}$} 
\EndFunction

\end{algorithmic}
\end{algorithm}

\subsubsection{Demonstration with example}

Consider the example \textsc{SumImpl} in Algorithm \ref{alg:sum8}, whose summation tree is illustrated in Figure \ref{fig:sum8}. We call Algorithm \ref{alg:optimized} with this \textsc{SumImpl} and $n=8$. First, the set of leaf nodes is $I=\{0,1,...,7\}$, where the smallest label is $i=0$. Next, the set $L_i=\{l^{i,j}: j\in I-\{i\}\}=\{2,4,4,6,6,8,8\}=\{2,4,6,8\}$ is computed. Examining the smallest value in $L_i$, we have $l=2$ and $J_l=\{j:j\in I-\{i\}\wedge l^{i,j}=l\}=\{1\}$. Therefore, \Call{BuildSubtree}{$\{1\}$} is called, reaching the stop condition. Then, the subtree with node \#0 and \#1 as its leaf nodes is constructed. The root of this subtree is denoted by $r$.

Next, examining the smallest value in the rest of $L_i$, we have $l=4$ and $J_l=\{j:j\in I-\{i\}\wedge l^{i,j}=l\}=\{2,3\}$. Therefore, \Call{BuildSubtree}{$\{2,3\}$} is called, where we have $I=\{2,3\}$, $i=2$, and $L_i=\{2\}$, and \Call{BuildSubtree}{$\{3\}$} is called there. \Call{BuildSubtree}{$\{2,3\}$} returns the subtree with node \#2 and \#3 as its leaf nodes. We then designate its root as the sibling node of $r$, and construct the parent node of this root and $r$. Then, the subtree with node \#0, \#1, \#2, and \#3 as its leaf nodes is constructed. $r$ is updated by the root of this subtree.

The next smallest value is $l=6$. We have $J_l=\{4,5\}$. Similarly, \Call{BuildSubtree}{$\{4,5\}$} is called, and it returns the subtree with node \#4 and \#5 as its leaf nodes. We merge its root with $r$, and construct the subtree with node \#0, \#1, \#2, \#3, \#4, and \#5 as its leaf nodes. $r$ is updated by the root of this subtree.

Finally, $l=8$ and $J_l=\{6,7\}$. \Call{BuildSubtree}{$\{6,7\}$} is called, and it returns the subtree with node \#6 and \#7 as its leaf nodes. We merge its root with $r$. Then the entire tree is constructed.

\subsubsection{Time complexity}
\label{sec:opt-time}

The time complexity of Algorithm \ref{alg:optimized} is $O(n^2t(n))$ and $\Omega(nt(n))$. The worst-case scenario occurs when adding $n$ summands in the right-to-left order. In this case, \Call{BuildSubtree}{} will be invoked with all suffixes of $\{0,1,...,n-1\}$, and $l^{i,j}$ for all $0\leq i<j<n$ will be calculated. The worst-case time complexity is $\Theta(n^2t(n))$. In practice, this order is cache-unfriendly, and thus no library uses it.

The best-case scenario corresponds to the sequential summation, where the summation tree will be constructed in one pass, and only  $l^{0,j}$ for all $0<j<n$ will be calculated. The best-case time complexity is $\Theta(nt(n))$. In practice, many libraries use similar orders, because these orders are cache-friendly and efficient.

\subsection{Adding support for matrix accelerators}

\subsubsection{Multi-term fused summation}

Matrix accelerators such as NVIDIA Tensor Cores are specialized hardware components in modern GPUs. Matrix accelerators enable assembly instructions that take a matrix $A=(a_{ij})_{M\times K}$, a matrix $B=(b_{ij})_{K \times N}$, and a matrix $C=(c_{ij})_{M\times N}$ as input, and produce a matrix $D=(d_{ij})_{M\times N}$ such that $D=A\times B+C$. The data types of $A$ and $B$ are identical. The data types of $C$ and $D$ are also identical, and their precision is no lower than the precision of $A$ and $B$. 

The numerical behavior of these assembly instructions remains undisclosed. Specifically, the computation of $d_{ij}=c_{ij}+\sum_{k=0}^{K-1}a_{ik}b_{kj}$ is executed in an undocumented way. Through delicate numerical experiments, prior works \cite{fasi_numerical_2021,li_fttn_2024} have found that for double-precision input, the computation is executed in a chain of standard FMAs; for low-precision input (specifically, when the precision of $A$ and $B$ is lower than float32), the computation of $d_{ij}=c_{ij}+\sum_{k=0}^{K-1}a_{ik}b_{kj}$ is executed based on multi-term fused summation:
\begin{itemize}
    \item The products are computed exactly, and the results are maintained in full precision without rounding after the multiplication.
    \item The summation of a group of summands is performed in a fixed-point manner. Specifically, the significands are aligned to the largest exponent of the summands, and then truncated to 24+ bits (i.e., no less than the precision of float32). The number of bits and the truncation method vary depending on the GPU architecture. 
    \item Then, the sum is converted to the floating-point number in the output data type of the instruction.
\end{itemize}

Note that the size of the group $w$, the width of the accumulator, and the detailed conversion method vary depending on the GPU architecture. In addition, previous works do not target the high-level APIs and libraries, so the accumulation orders of them remains unknown.

Our proposed solutions can work for standard FMAs. However, multi-term fused summation requires a new method, because it is executed in a non-standard, IEEE-754-incompliant way. Specifically, in multi-term fused summation, $w$ summands (e.g., $x_0=c$, and $x_i=a_{i-1} b_{i-1}$ for $1\le i<w$) are summed in a fixed-point manner, thus making the result independent of the summation order. To represent this operation in the summation tree, we should use a node with $w$ children instead of a node with two children. Therefore, the summation tree should be an $w$-way tree.


\subsubsection{Constructing the multiway summation tree}

To adapt to the multiway tree, we first revisit BasicFPRev in Section \ref{sec:method}. The first two steps still work because we find that the key equation $l^{i,j}=n-\Call{SumImpl}{A^{i,j}}$ remains valid in multi-term fused summation. Thus, the values of $l^{i,j}$ can be obtained in the same way, and we only need to redesign the tree construction algorithm in the third step.

Then, we revisit the tree construction algorithm in Algorithm \ref{alg:optimized}. In \Call{BuildSubtree}{$I$}, we calculate $l^{i,j}$ for a fixed $i$ and all $j\in I-\{i\}$, enumerate them in ascending order, and maintain $r$ as the root of the largest constructed subtree containing node \#$i$. For some $l\in L_i=\{l^{i,j}:j\in I-\{i\}\}$ and $J_l=\{j : j\in I-\{i\} \land l^{i,j}=l\}$, the return value of \Call{BuildSubtree}{$J_l$} is the subtree with $J_l$ as its leaf nodes. The root of this subtree must be the sibling node of $r$, so we can create a new node as their parent node and update $r$. However, this relation is not always true for the multiway tree.

In addition to being sibling node, the root of the subtree may also be the parent node of $r$ in the multiway tree. For example, suppose a 5-way tree with leaf nodes $I=\{0,1,2,3,4\}$ as the children of the root. Then, when $r=0$, $l=5$, and $J_l=\{1,2,3,4\}$, solving the subproblem for $J_l$ should return a partial subtree with $J_l$ as its leaves. The root node of the subtree is the parent node of $r$.

To distinguish the two cases, we observe the return value of \Call{BuildSubtree}{$J_l$}, denoted by $T'$, and the complete subtree rooted at the root of $T'$, denoted by $T_c$. In the first case, the root of $T'$ should be the sibling of $r$, and $T'=T_c$. In the second case, the root of $T'$ should be the parent of $r$, and $T' \subset T_c$. Therefore, we can compare the size of $T'$ (denoted by $n_\mathrm{leaves}^{T'}$) with the size of $T_c$ (denoted by $n_\mathrm{leaves}^{T_c}$). We note that $n_\mathrm{leaves}^{T'}=|J_l|$, and $n_\mathrm{leaves}^{T_c}=\max \{l^{j,k} : j,k \in J_l\}=\max(L_{\min(J_l)})$. Therefore, if $\max(L_{\min(J_l)})=|J_l|$, then the root of $T'$ should be the sibling of $r$, so we should create a new node as their parent node and update $r$ with the index of this new node. Otherwise, $\max(L_{\min(J_l)})>|J_l|$, so the root of $T'$ should be the parent of $r$, and thus we should add an edge from $r$ to the root of $T'$, and update $r$ with the root.

Through this modification, the multiway tree can be correctly constructed. We elaborate on the above process in Algorithm \ref{alg:fprev}, i.e., the algorithm of FPRev. It has the same time complexity as Algorithm \ref{alg:optimized} (note that Algorithm \ref{alg:optimized} just corresponds to the case where $\max(L_{\min(J_l)})=|J_l|$), and supports multi-term fused summation. 

\begin{algorithm}[ht]
\caption{The algorithm of FPRev.}
\label{alg:fprev}
\begin{algorithmic}
\Require Implementation \textsc{SumImpl} and the number of summands $n$
\Ensure Summation tree of \textsc{SumImpl}

\Function{FPRev}{\textsc{SumImpl},$n$}
    \Function{BuildSubtree}{$I$}
    \State $T\gets \varnothing$
    \If{$|I|=1$} \Comment{stop condition}
        \State \Return $(T,1)$
    \EndIf
    \State $i \gets \min(I)$; $L_i\gets \varnothing$
    \For{$j\in I-\{i\}$} \Comment{calculate $l^{i,j}$ on demand}
        \State $A^{i,j} \gets (1,1,...,1)$; $A^{i,j}_i\gets M$; $A^{i,j}_j\gets -M$ 
        \State $l^{i,j}\gets n-\Call{SumImpl}{A^{i,j}}$
        \State $L_i\gets L_i\cup \{l^{i,j}\}$
    \EndFor
    \State $r \gets i$ \Comment{current root of the subtree}
    \For{$l \in L_i$ in ascending order}
        \State $J_l \gets \{j : j\in I-\{i\} \land l^{i,j}=l\}$
        \State $(T',n_\mathrm{leaves}^{T_c}) \gets \Call{BuildSubtree}{J}$ 
        \State $T\gets T\cup T'$
        \If{$|J_l|=n_\mathrm{leaves}^{T_c}$} \Comment{$T'=T_c$}
            \State $T\gets T\cup\{(r,r+n), (\mathrm{GetRoot}(T'),r+n)\}$
            \State $r \gets r+n$
        \Else \Comment{$T' \subset T_c$}
            \State $T\gets T\cup\{(r,\mathrm{GetRoot}(T'))\}$
            \State $r \gets \mathrm{GetRoot}(T')$
        \EndIf
    \EndFor
    \State \Return $(T,\max(L_i))$
    \EndFunction
\State $(T,n_\mathrm{leaves}^{T_c})\gets \Call{BuildSubtree}{\{0,1,...,n-1\}}$ 
\State \Return $T$
\EndFunction

\end{algorithmic}
\end{algorithm}

\subsection{Time complexity and correctness}

Following the same analysis in Section \ref{sec:opt-time}, the time complexity of FPRev is $O(n^2t(n))$ and $\Omega(nt(n))$, and the probability of the worst-case time complexity $O(n^2t(n))$ is low. The correctness of it is also guaranteed by design and can be proven following the same process in Section \ref{sec:basic-timecor}. 
\section{Case study}
\label{sec:casestudy}

In this section, we apply FPRev to two prevalent numerical libraries: NumPy \cite{harris_array_2020}, the most popular Python library for numerical computing on CPUs, and PyTorch \cite{paszke_pytorch_2019}, a very popular Python library for numerical computing on GPUs. We successfully identify and analyze the undocumented accumulation orders in these libraries.

\subsection{NumPy's implementation on CPUs}
\label{sec:casenumpy}

We use FPRev to test NumPy (version 1.26) on three CPUs:\begin{itemize}
    \item CPU-1: Intel Xeon E5-2690 v4 (24 v-cores)
    \item CPU-2: AMD EPYC 7V13 (24 v-cores)
    \item CPU-3: Intel Xeon Silver 4210 (40 v-cores)
\end{itemize}


\textbf{Summation.} On these CPUs, we find that NumPy implements identical accumulation order for the summation function in single precision. Therefore, Numpy's summation is verified to be reproducible across these systems, and can be used in software requiring numerical reproducibility.

The accumulation order is sequential for $n<8$. For $8\leq n\leq 128$, NumPy implements an eight-way summation. Each way $i$ sums up $a_i, a_{i+8}, a_{i+16}, ...$ sequentially, and the sums of eight ways are summed using pairwise summation. For example, Figure \ref{fig:numpy_sum_32} shows the accumulation order for $n=32$. This accumulation order implies that developers can leverage the eight-way SIMD instructions in the CPU to accelerate computation. For $n>128$, NumPy increases the number of ways, thus leveraging multi-threading for large-scale summation.

\begin{figure}[h]
    \centering
    \begin{subfigure}[b]{0.6\linewidth}
    \includegraphics[width=1\linewidth]{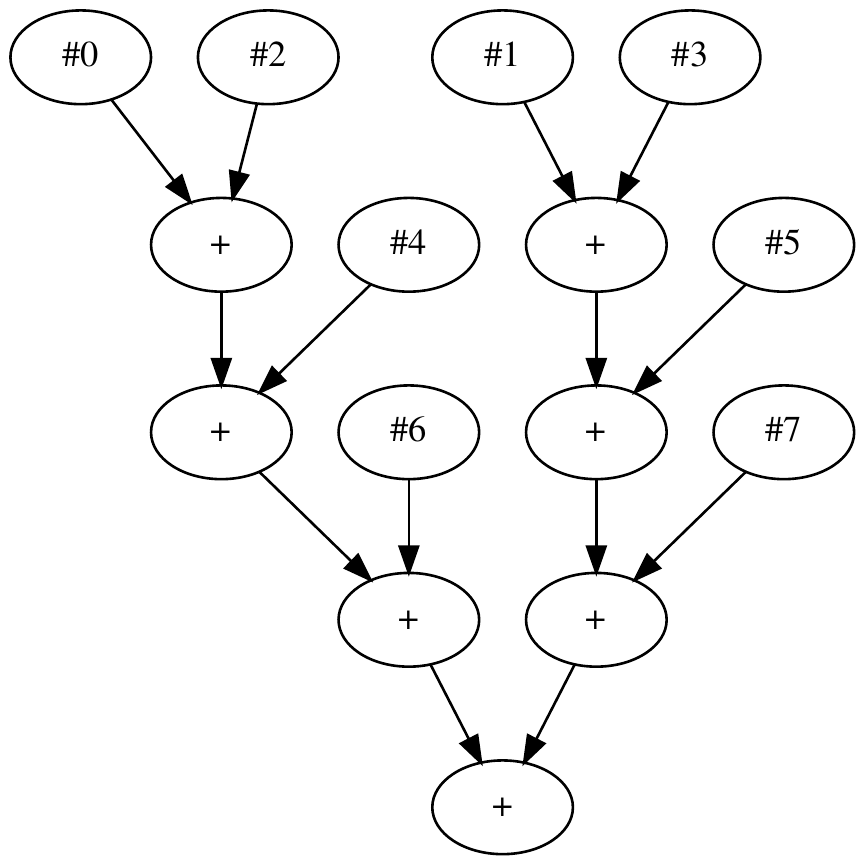}
    \caption{On Intel Xeon E5-2690 v4 and AMD EPYC 7V13 (24 v-cores).}
    \end{subfigure}
    \begin{subfigure}[b]{0.39\linewidth}
    \includegraphics[width=1\linewidth]{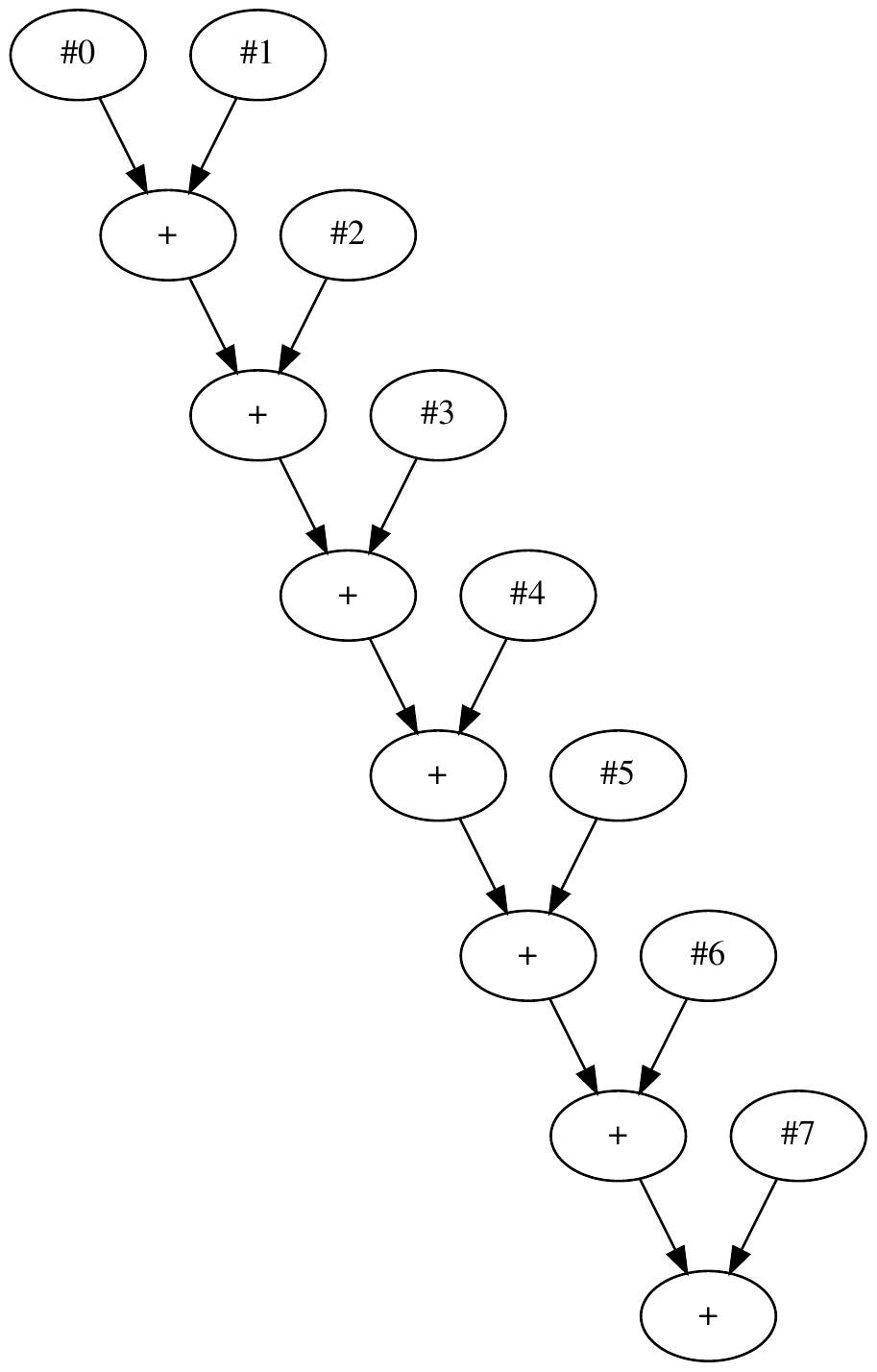}
    \caption{On Intel Xeon Silver 4210 (40 v-cores).}
    \end{subfigure}

\caption{The accumulation orders of NumPy's $8\times 8$ matrix-vector multiplication on different CPUs.}
\label{fig:numpy_gemv}
\end{figure}

\textbf{Other AccumOps.} We also test NumPy's dot product, matrix-vector multiplication, and matrix multiplication functions in single precision. We observe discrepancies in the accumulation order across the tested CPUs. For example, Figure \ref{fig:numpy_gemv} shows the accumulation orders of Numpy's $n\times n$ matrix-vector multiplication for $n=8$ on the CPUs. We note that on CPU-1 and CPU-2, the 32 products of each output element are accumulated using 2-way summation, whereas on CPU-3, which has more cores than CPU-1 and CPU-2, the products are accumulated sequentially. 

In summary, NumPy's summation function is safe for developing reproducible software for these CPUs, while other AccumOps of NumPy should not be used in software requiring numerical reproducibility.

\subsection{PyTorch's implementation on GPUs}
\label{sec:casepytorch}

We use FPRev to test PyTorch (version 2.3) on three GPUs:
\begin{itemize}
    \item GPU-1: NVIDIA V100 (5120 CUDA cores)
    \item GPU-2: NVIDIA A100 (6912 CUDA cores)
    \item GPU-3: NVIDIA H100 (16896 CUDA cores)
\end{itemize}

On these GPUs, we observe findings similar to those for NumPy: PyTorch implements identical accumulation orders for the summation function in single precision but not for the BLAS operations. Therefore, PyTorch's summation function is safe for developing reproducible software for these GPUs, while other AccumOps of PyTorch should not be used in software requiring numerical reproducibility.



\begin{figure}[h]
    \centering
    \begin{subfigure}[b]{\linewidth}
    \includegraphics[width=1\linewidth]{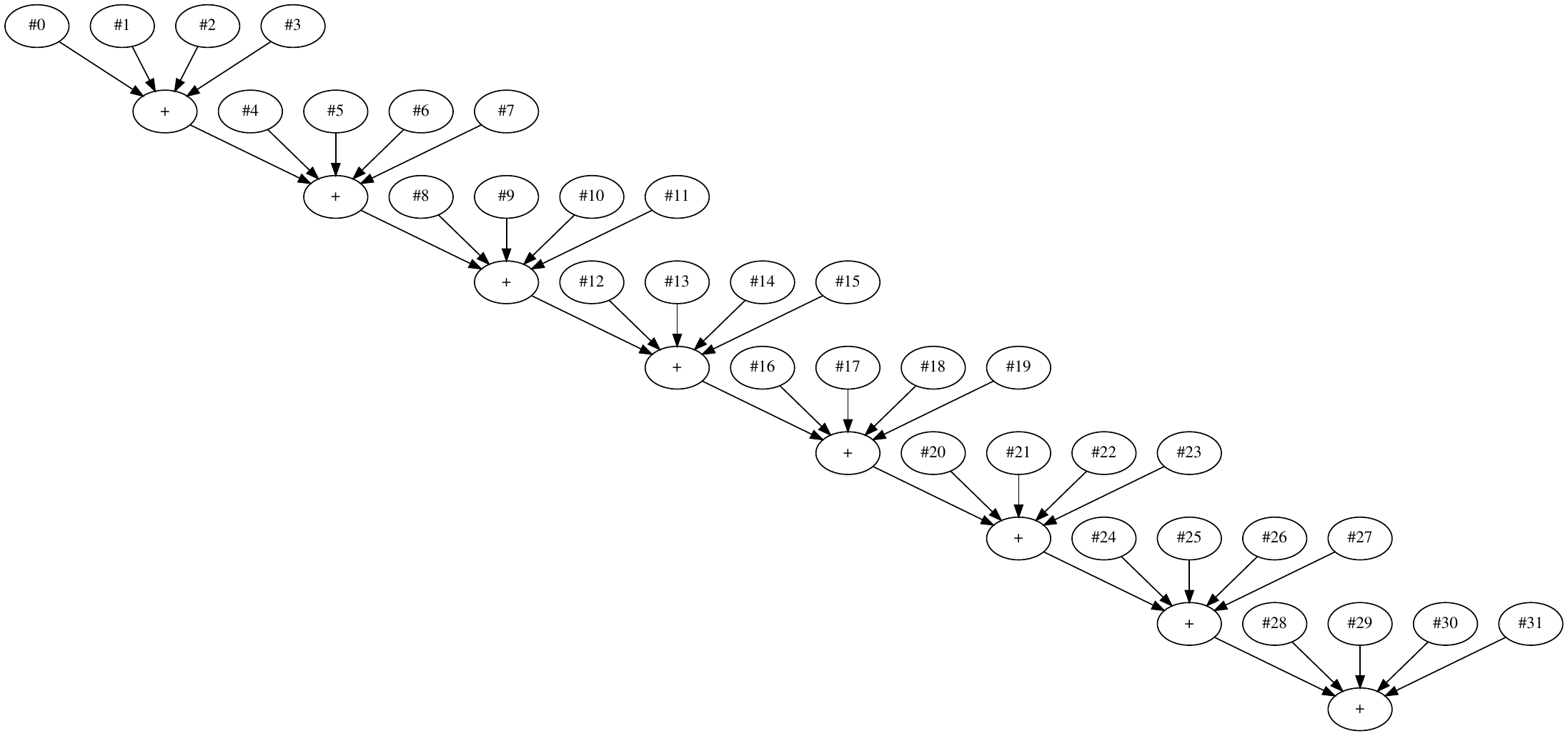}
    \caption{On NVIDIA V100.}
    \end{subfigure}
    \begin{subfigure}[b]{\linewidth}
    \includegraphics[width=1\linewidth]{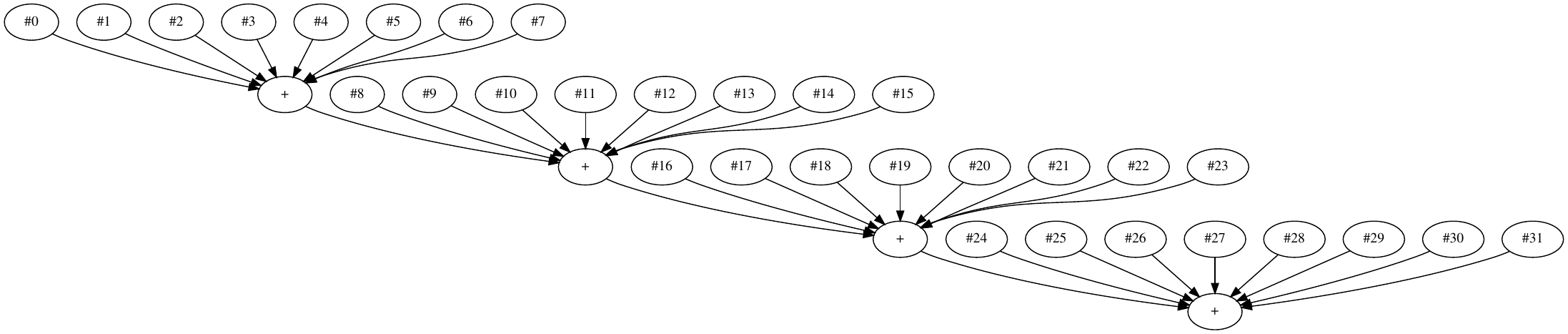}
    \caption{On NVIDIA A100.}
    \end{subfigure}

    \begin{subfigure}[b]{\linewidth}
    \includegraphics[width=1\linewidth]{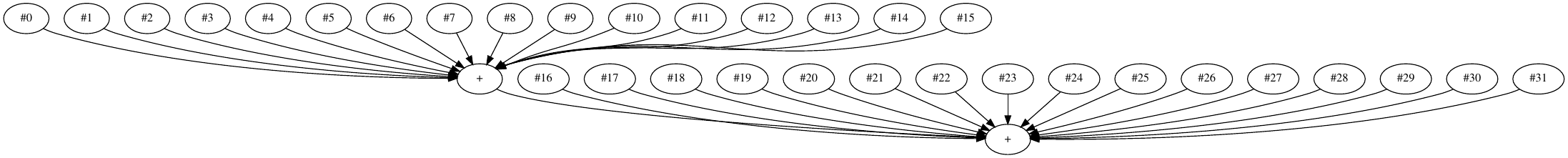}
    \caption{On NVIDIA H100.}
    \end{subfigure}
\caption{The accumulation orders of PyTorch's half-precision $32\times 32\times 32$ matrix multiplication on Tensor Cores.}
\label{fig:tc_gemm}
\end{figure}

\textbf{Matrix multiplication on Tensor Cores.} To enable Tensor Core computation, we apply FPRev to half-precision matrix multiplication in PyTorch, which is implemented using the cuBLAS backend. The results show that the summation tree is a 5-way tree on NVIDIA V100, a 9-way tree on A100, and a 17-way tree on H100, corroborating the conclusion in \cite{fasi_numerical_2021,li_fttn_2024}, which states that the Tensor Cores on NVIDIA Volta, Ampere, and Hopper architectures use (4+1)-, (8+1)-, and (16+1)-term fused summation respectively. For example, Figure \ref{fig:tc_gemm} shows the summation trees for $n=32$ on these devices. 

We also examine the SASS assembly instructions they use, and observe that V100 uses the \texttt{HMMA.884} instruction, and both A100 and H100 use the \texttt{HMMA.16816} instruction. \textbf{Interestingly, an \texttt{HMMA.16816} instruction on A100 indicates the shape of the inputs where the accumulation dimension $K=16$, but it is implemented through (8+1)-term fused summation by the A100 Tensor Core hardware.}

\section{Performance evaluation}
\label{sec:evaluation}

\subsection{Experiment design}

\begin{figure*}[t]
    \centering
    \includegraphics[width=\linewidth]{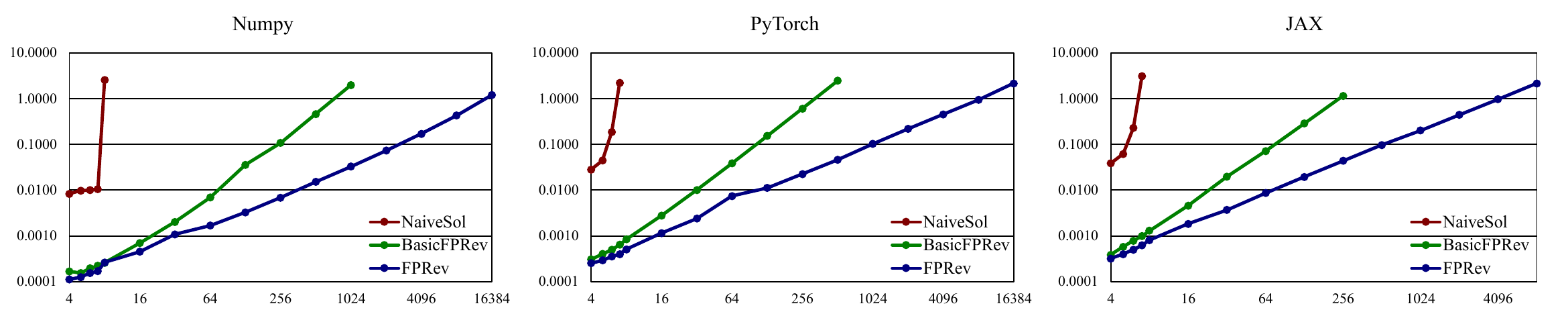}
    \caption{Execution time of applying NaiveSol, BasicFPRev, and FPRev to the summation functions in NumPy, PyTorch, and JAX. The vertical axis represents execution time in seconds. The horizontal axis represents the number of summands $n$.}
    \label{fig:rq1}
\end{figure*}
\begin{figure*}[t]
    \centering
    \includegraphics[width=\linewidth]{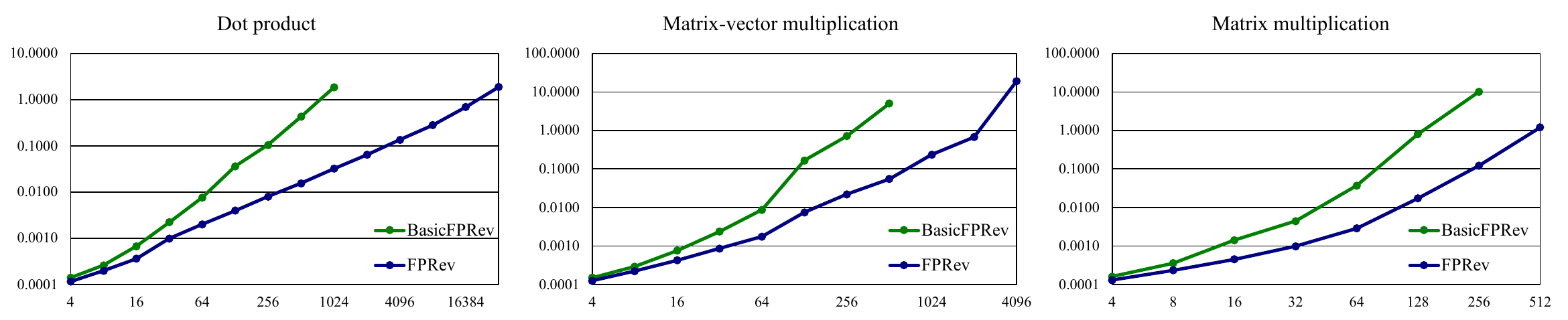}
    \caption{Execution time of applying BasicFPRev and FPRev to the dor product, matrix-vector multiplication, and matrix multiplication functions in NumPy. The vertical axis represents execution time in seconds. The horizontal axis represents the number of summands $n$.}
    \label{fig:rq2}
\end{figure*}

In this section, we evaluate the efficiency of FPRev. Specifically, we aim to answer the following research questions (RQs): 
\begin{itemize}
    \item RQ1: how efficient is FPRev when applied to different libraries?
    \item RQ2: how efficient is FPRev when applied to different operations?
    \item RQ3: how efficient is FPRev on different CPUs and GPUs?
\end{itemize}

To answer the RQs, we measure the execution time (wall-clock time) of applying our solutions to tested libraries and operations. We implement FPRev (Algorithm \ref{alg:fprev}) in Python (version 3.11). For comparison, we also implement the basic solution (Algorithm \ref{alg:basic}), denoted by BasicFPRev.

\subsection{RQ1: How efficient is FPRev when applied to different libraries?}
\label{sec:rq1}

\begin{figure*}[t]
    \centering
    \includegraphics[width=\linewidth]{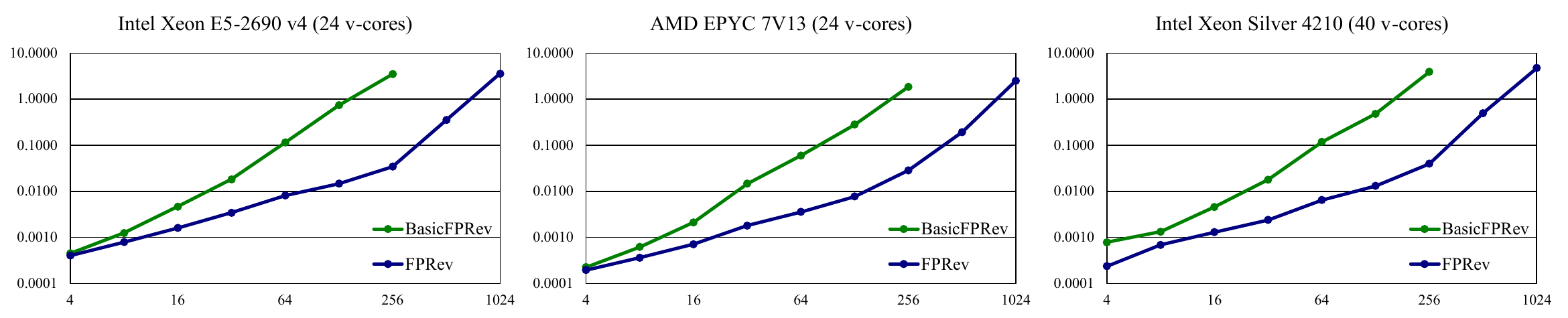}
    \includegraphics[width=\linewidth]{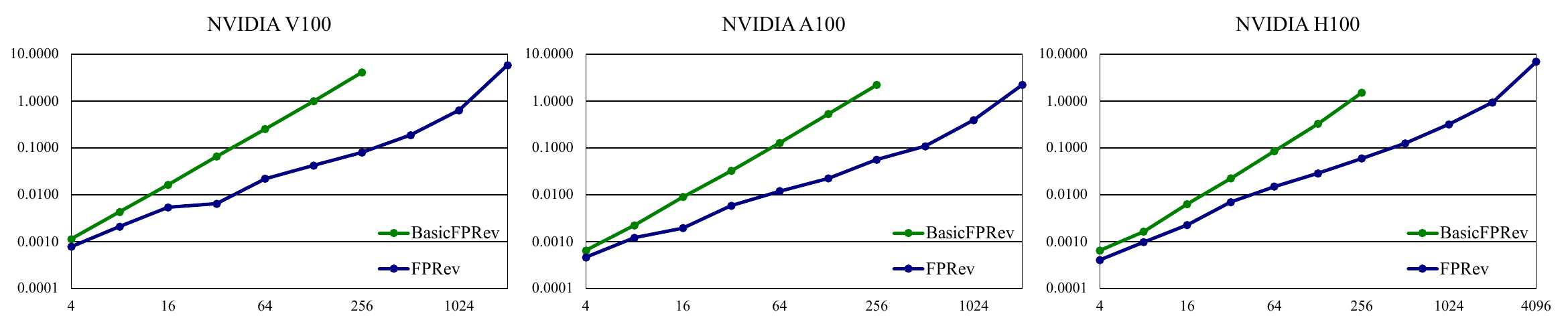}
    \caption{Execution time of applying BasicFPRev and FPRev to the matrix multiplication functions in PyTorch on different CPUs and GPUs. The vertical axis represents execution time in seconds. The horizontal axis represents the number of summands $n$.}
    \label{fig:rq3}
\end{figure*}

For RQ1, we test the single-precision summation function in three libraries: NumPy (version 1.26) \cite{harris_array_2020}, PyTorch (version 2.3) \cite{paszke_pytorch_2019}, and JAX (version 0.4) \cite{bradbury_jax_2018}. We run the experiments on Intel Xeon E5-2690 v4 with 24 v-cores. In these experiments, we also implement the naive brute-force solution (Section \ref{sec:naive}), denoted by ``NaiveSol'', to show its extremely low efficiency. For remaining RQs, we omit the naive solution because it is proven to be too inefficient.

We begin with the number of summands $n=4$, and increment $n$ until the execution time exceeds one second. Each experiment is carried out 10 times, and the arithmetic mean of the 10 results is reported in Figure \ref{fig:rq1}. The red curves indicate that the execution time of NaiveSol grows exponentially as $n$ grows. The results substantiate the $O(4^n/n^{3/2}\cdot t(n))$ time complexity of NaiveSol. The green and blue lines show that the execution time of BasicFPRev and FPRev grows polynomially. The different slopes also demonstrate that the execution time of BasicFPRev is longer than that of FPRev, and increases more rapidly as $n$ increases. This is because the time complexity of BasicFPRev is $\Theta(n^2t(n))$, while that of FPRev is $\Omega(nt(n))$ and $O(n^2t(n))$.

These trends suggest that the scalability of BasicFPRev is much better than that of NaiveSol, and the scalability of FPRev is even better. For example, when $n=16$, NaiveSol can take over 24 hours to produce an output, but BasicFPRev and FPRev only take less than 0.01 seconds. If $n=8192$, BasicFPRev will take over 100 seconds to produce an output, but FPRev only takes about 1 second.




\subsection{RQ2: How efficient is FPRev when applied to different operations? }
\label{sec:rq2}

For RQ2, we test the single-precision dot product, matrix-vector multiplication, and matrix multiplication functions in NumPy on Intel Xeon E5-2690 v4. The time complexity of these functions is $O(n)$, $O(n^2)$, and $O(n^3)$, respectively.

Similarly, we begin with $n=4$, and increment $n$ until the execution time exceeds one second. Each experiment is carried out 10 times, and the arithmetic mean of the 10 results is reported in Table \ref{fig:rq2}. The different slopes indicate that the time complexity of BasicFPRev is higher than that of FPRev. In addition, as the time complexity of the workload increases, the growth speed of the runtime with regard to $n$ accelerates. Therefore, the speedup of FPRev over BasicFPRev is more pronounced as the workload is more complex. For example, for $n=256$, FPRev is $13.0\times$ as fast as BasicFPRev for dot product, $32.3\times$ for matrix-vector multiplication, and $82.1\times$ for matrix multiplication.


\subsection{RQ3: How efficient is FPRev on different CPUs and GPUs?}
\label{sec:rq3}

For RQ3, we test the single-precision matrix multiplication in PyTorch on the CPUs and GPUs listed in Section \ref{sec:casestudy}. Similarly, we begin with $n=4$, and increment $n$ until the execution time exceeds one second. Each experiment is carried out 10 times, and the arithmetic mean of the 10 results is reported in Figure \ref{fig:rq3}. The results demonstrate consistent improvements in the runtime of FPRev compared to BasicFPRev. Therefore, FPRev is consistently more efficient than BasicFPRev on different devices.

\section{Discussion}

\subsection{Limitation and mitigation}

\subsubsection{Dynamic range of the input data type} 

When applying FPRev to data types with low dynamic range, the mask $M$ may be too small to effectively mask the sum of ones. For example, the maximum value of the 8-bit floating-point number with 4 exponent bits (FP8-e4m3) defined in \cite{micikevicius_fp8_2022} is $1.75 \times 2^8$, so the condition $M\gg n$ may not hold, and $\pm M+\sigma \neq \pm M$ for $0\leq \sigma \leq n-2$. To mitigate the issue, we can replace the ones in the masked all-one arrays with smaller numbers (e.g., $2^{-9}\times2^{-9}$ for FP8-e4m3 matrix multiplication), and scale the sum back to an integer between 0 and $n-2$ when calculating $l^{i,j}$. This solution does not affect efficiency.

\subsubsection{Precision of the accumulator} 

The precision of the floating-point accumulator can limit the input size that FPRev supports. For example, float32 has a precision of 24 bits, so the maximum number of summands ($n$) that FPRev supports is $2^{24}+1=16777217$ for float32 accumulation operations. For larger numbers, the sum of $n-2$ ones cannot be represented precisely in float32, so the sum of masked all-one arrays may be incorrect. This issue can be mitigated by dynamically replacing the multiple ones corresponding to a constructed subtree with one and multiple zeros, as if compressing the constructed subtree into one node.

Specifically in the \Call{BuildSubtree}{$I$} function of FPRev (Algorithm \ref{alg:fprev}), the computation of $\textsc{sum}(A^{i,j})=0$ is accurate for $j$s such that $l^{i,j}=n$, so we can build the subtree for them in the end. We extract those $j$s to $J = \{j :  l^{i,j}=n\}$ after computing $L_i = \{l^{i,j} : j \in I-\{i\}\}$. Then, we set the values at $J$ to 0, and build the subtree for $I-J$ recursively (which results in a smaller subproblem). After the subtree for $I-J$ is constructed, we set the values at $J$ back to 1.0, and set the values at $I-J-\{i\}$ to 0. Now, the constructed subtree (containing $I-J$) is treated as a node, represented by node $\#i$. Next, we run the original tree construction algorithm (the last iteration of the for-loop, i.e., $l=|All|$) for $J$, and then the whole tree is constructed.

Combining the two mitigation techniques, the modified version of FPRev is shown in Algorithm \ref{alg:modified}. This version is applicable to data types with low dynamic range and low accumulation precision, such as 16-bit and 8-bit floating-point formats (including BF16, FP16, FP8-e5m2, and FP8-e4m3 on recent Tensor Cores).

\begin{algorithm}[ht]
\caption{The modified version of FPRev.}
\label{alg:modified}
\begin{algorithmic}
\Require Implementation \textsc{SumImpl}, the number of summands $n$, large value $M$, and tiny value $e$
\Ensure Summation tree of \textsc{SumImpl}

\Function{ModifiedFPRev}{$\textsc{SumImpl},n, M, e$}
    \Function{BuildSubtree}{$I, All$}
    \State $T\gets \varnothing$
    \If{$|I|=1$} \Comment{stop condition}
        \State \Return $T$
    \EndIf
    \State $i \gets \min(I)$; $L_i\gets \varnothing$
    \For{$j\in I-\{i\}$} \Comment{calculate $l^{i,j}$ on demand}
        \State $A^{i,j}_k \gets e$ for $k\in All$; $A^{i,j}_i\gets M$; $A^{i,j}_j\gets -M$ 
        \State $l^{i,j}\gets |All|-\Call{SumImpl}{A^{i,j}}/e$
        \State $L_i\gets L_i\cup \{l^{i,j}\}$
    \EndFor
    \State $J\gets \{j:l^{i,j}=\max(L_i)\}$
    \State $A^{i,j}_k \gets 0$ for $k\in J$ \Comment{Ignoring $J$}
    \State $(T, n_\mathrm{leaves}^{T_{c}}) \gets \textsc{BuildSubtree}(I-J, All-J)$
    \State $r \gets \mathrm{GetRoot}(T)$ 
    \State $A^{i,j}_k \gets e$ for $k\in J$
    \State $K\gets I-J-\{i\}$
    \State $A^{i,j}_k \gets 0$ for $k\in K$ \Comment{Treating $I-J$ as $\{i\}$}
    \State $(T',n_\mathrm{leaves}^{T_c}) \gets \Call{BuildSubtree}{J, All-K}$ 
    \State $T\gets T\cup T'$
    \If{$|J|=n_\mathrm{leaves}^{T_c}$} \Comment{$T'=T_c$}
        \State $T\gets T\cup\{(r,r+n), (\mathrm{GetRoot}(T'),r+n)\}$
    \Else \Comment{$T' \subset T_c$}
        \State $T\gets T\cup\{(r,\mathrm{GetRoot}(T'))\}$
    \EndIf
    \State \Return $(T,\max(L_i))$
    \EndFunction
\State $All\gets \{0,1,...,n-1\}$
\State $(T,n_\mathrm{leaves}^{T_c})\gets \Call{BuildSubtree}{All,All}$ 
\State \Return $T$
\EndFunction

\end{algorithmic}
\end{algorithm}



\subsection{Extensibility and future work}
\label{sec:scope}

Section \ref{sec:casestudy} has demonstrated that FPRev can be applied to popular numerical libraries like NumPy and PyTorch. FPRev can be applied to other accumulation implementations as long as they fall within the scope detailed in Section \ref{sec:problem}. In practice, we find most popular libraries have deterministic reduction orders and fall into the scope. FPRev also works for accumulation operations in collective communication primitives, such as the AllReduce operation, if their accumulation order is predetermined.

To further extend our tool to other functions based on special summation algorithms, the algorithms must satisfy the property $l^{i,j}=n-\Call{SumImpl}{A^{i,j}}$ or its variant in Algorithm \ref{alg:modified}. For example, the next generation of Tensor Core will support the microscaling data format \cite{rouhani_microscaling_2023}, including the 4-bit and 6-bit formats MXFP4 and MXFP6. If their dynamic range and accumulator precision permit and the property holds, our methods can reveal the accumulation order within a block of microscaling numbers. Then, we can treat a block as one summand, and use FPRev to construct the summation tree for the summation of the blocks, and then expand each block to a subtree. 

In addition to revealing accumulation orders, we plan to extend our methods to detect more floating-point behaviors in matrix accelerators. For example, we can determine the rounding mode and the precision of the accumulator of Tensor Cores by enumerating $n=1,2,...$ and checking the result of $2^n+1.75-2^n$. We are designing more numerical experiments to identify how the block fused multiply-add is conducted. Then, with the information detected, we can model the exact behavior of the hardware matrix accelerators.

Another direction is further optimizing the efficiency of FPRev. For example, we can randomize the selection of $i\in I$ in the FPRev algorithm, as if selecting the random pivot in quick sort. This might reduce the expected time complexity of FPRev.
\section{Conclusion}
In this paper, we introduce FPRev, a diagnostic tool for revealing the accumulation order in software and hardware implementations through numerical testing. It can help verify and facilitate the development of reproducible software. As a case study, FPRev reveal the undisclosed accumulation orders of prevalent numerical libraries such as NumPy and PyTorch on different CPUs and GPUs. We also demonstrated the efficiency FPRev through experiments covering various implementations and devices. Our source code is available at \url{https://github.com/peichenxie/FPRev}, encouraging further investigation and improvement by the research community.


\bibliographystyle{plainurl}
\bibliography{references}

\appendix
\section{Artifact Appendix}

\subsection*{Abstract}

The repository includes the source code of FPRev and the source code for reproducing the experiments of the paper. The following content includes main claims that can be verified via experiments, and instructions to reproduce the experiments.

\subsection*{Scope}

The main claims include:

\begin{enumerate}
    \item FPRev is functional to reveal the floating-point accumulation orders in common implementations. This claim is detailed by Section \ref{sec:casestudy} ``Case study'' in the paper. To verify this claim, run \texttt{python experiments/casestudy.py} and check the output files in the \texttt{outputs} directory.
    \begin{enumerate}
        \item \texttt{outputs/Numpy*.pdf} represents the revealed accumulation orders for NumPy, as discussed in Section \ref{sec:casenumpy} ``NumPy's implementation on CPUs''. Among them, \texttt{outputs/NumpyGEMV8.pdf} corresponds to Figure \ref{fig:numpy_gemv} of the paper, if the CPU models are consistent to those in the paper.
        \item \texttt{outputs/Torch*.pdf} represents the revealed accumulation orders for PyTorch, as discussed in Section \ref{sec:casepytorch} ``PyTorch's implementation on GPUs''. Among them, \texttt{outputs/TorchF16GEMM32.pdf} corresponds to Figure \ref{fig:tc_gemm} of the paper, if the GPU models are consistent to those in the paper.
    \end{enumerate}
    \item FPRev is efficient. This claim is detailed by Section \ref{sec:evaluation} ``Performance evaluation'' in the paper. To verify this claim, run \texttt{python experiments/rq1.py}, \texttt{python experiments/rq2.py}, and \texttt{python experiments/rq3.py}, and check the output files in the \texttt{outputs} directory.
    \begin{enumerate}
        \item \texttt{outputs/rq1.csv} provides the results of Section \ref{sec:rq1} ``RQ1: How efficient is FPRev when applied to different libraries?''. It corresponds to Figure \ref{fig:rq1} if the CPU model is consistenst to that in the paper.
        \item \texttt{outputs/rq2.csv} provides the results of Section \ref{sec:rq2} ``RQ2: How efficient is FPRev when applied to different operations?''. It corresponds to Figure \ref{fig:rq2} if the CPU model is consistenst to that in the paper.
        \item \texttt{outputs/rq3.csv} provides the results of Section \ref{sec:rq3} ``RQ3: How efficient is FPRev on different CPUs and GPUs?''. It corresponds to Figure \ref{fig:rq3} if the CPU and GPU models are consistenst to those in the paper.
    \end{enumerate}
\end{enumerate}

\subsection*{Contents}

\paragraph{Installation}

\begin{verbatim}
sudo apt install graphviz
git clone https://github.com/peichenxie/FPRev.git
cd FPRev
pip install .
pip install -r experiments/requirements.txt
\end{verbatim}

\paragraph{Running experiments}

\begin{itemize}
    \item To reproduce the results in Section \ref{sec:casestudy} (Case study), run \texttt{python experiments/casestudy.py} on different hardware models.
    \item To reproduce the results in Section \ref{sec:rq1} (RQ1: How efficient is FPRev when applied to different libraries?), run \texttt{python experiments/rq1.py}.
    \item To reproduce the results in Section \ref{sec:rq2} (RQ2: How efficient is FPRev when applied to different operations?), \texttt{run python experiments/rq2.py}.
    \item To reproduce the results in Section \ref{sec:rq3} (RQ3: How efficient is FPRev on different CPUs and GPUs?), run \texttt{python experiments/rq3.py} on different hardware models.
\end{itemize}

Then, check the output files in the outputs directory. See \texttt{outputs/README.md} for more information.

\subsection*{Hosting}

\url{https://github.com/peichenxie/FPRev}

\subsection*{Requirements}

The artifact requires the following platform:

\begin{itemize}
    \item GPU: NVIDIA V100 or newer
    \item OS: Ubuntu 22.04
    \item Software: Python (version 3.11)
\end{itemize}

If you wish to reproduce the results in the paper, please use the identical CPU and GPU models:
\begin{enumerate}
    \item CPU: Intel Xeon E5-2690 v4 (24 v-cores), GPU: NVIDIA V100 (5120 CUDA cores)
    \item CPU: AMD EPYC 7V13 (24 v-cores), GPU: NVIDIA A100 (6912 CUDA cores)
    \item CPU: Intel Xeon Silver 4210 (40 v-cores), GPU: NVIDIA H100 (16896 CUDA cores)
\end{enumerate}

\end{document}